\documentclass{amsart}
\usepackage{amsmath}
\usepackage{graphicx}
\usepackage{amssymb}
\usepackage[margin=1in]{geometry}
\usepackage{url}
\usepackage{xcolor}
\usepackage{hyperref}

\newtheorem{thm}{Theorem}[section]

\theoremstyle{definition}

\theoremstyle{remark}
\newtheorem{rem}[thm]{Remark}

\begin{document}

\title[Kinetic Aleksandrov estimate]{An Aleksandrov estimate for kinetic equations with rough coefficients}
\author[A. Loher]{Am\'{e}lie Loher}
\address{All Souls College, University of Oxford}
\email{\tt amelie.loher@all-souls.ox.ac.uk}
\author[C. Mooney]{Connor Mooney}
\address{Department of Mathematics, UC Irvine}
\email{\tt mooneycr@uci.edu}
\author[C. Mouhot]{Cl\'{e}ment Mouhot}
\address{DPMMS, University of Cambridge}
\email{\tt c.mouhot@dpmms.cam.ac.uk}
\author[L. Niebel]{Lukas Niebel}
\address{Department of Mathematics, ETH Z\"{u}rich}
\email{\tt lukas.niebel@math.ethz.ch}
\author[L. Silvestre]{Luis Silvestre}
\address{Mathematics Department, University of Chicago}
\email{luiss@uchicago.edu}
\author[Y. Zhu]{Yuzhe Zhu}
\address{SSE, The Chinese University of Hong Kong, Shenzhen}
\email{zhuyuzhe@cuhk.edu.cn}
\begin{abstract}
  We give an exposition of a proof, discovered by Astra, of an Aleksandrov estimate for kinetic equations with rough coefficients. We explain this method in the simpler elliptic and parabolic settings as well.
\end{abstract}
\maketitle

\section{Introduction}

The Aleksandrov–Bakelman–Pucci estimate, also called the Aleksandrov maximum principle in the Soviet school, plays a central role in the regularity theory for elliptic and parabolic equations in nondivergence form. Its proof starts from a geometric insight by Aleksandrov in the context of the Monge--Amp\`ere equation \cite{MR0147776}, relating the oscillation of a convex function with the set of its supporting planes. It was applied to the study of elliptic equations in nondivergence form by Bakel'man \cite{MR126604} and Pucci \cite{MR214905,MR208150}. Krylov later extended this estimate to the parabolic setting \cite{MR420016}.

The proof of Hölder continuity and the Harnack inequality by Krylov and Safonov~\cite{KS} is based on three main nontrivial ingredients:
\begin{itemize}
    \item[(A)] An $\epsilon$-density argument, which is usually deduced from the ABP estimate in the non-divergence-form theory. This is equivalent to what is called \emph{the first growth lemma} in the school of Landis \cite{Landis}. This estimate is a ``weak'' measure-to-pointwise result that says, e.g.\ for supersolutions, that being large enough in \emph{most} of a cylinder implies pointwise positivity in a smaller future cylinder. At its core this step is a universal non-concentration result. It corresponds to the first lemma of De Giorgi in the divergence-form theory.
    \item[(B)] An expansion of positivity lemma, which is obtained by constructing a barrier function. This is sometimes called \emph{the double-ball property}.
    \item[(C)] A covering argument, which originally took the form of the growing ink-spots lemma. In some proofs, it is replaced with the Calder\'on--Zygmund decomposition.
\end{itemize}
The ingredients (B) and (C) together allow us to improve the $\epsilon$-density argument (weak measure-to-pointwise estimate) into a \emph{strong} measure-to-pointwise estimate. The latter says, e.g. for supersolutions, that being large enough in \emph{any fraction} of a cylinder implies pointwise positivity in a smaller future cylinder. This is a substantial upgrade of
(A) and really captures a degree of smoothing: it readily implies a reduction of oscillation and Hölder continuity. This result is called the \emph{second growth lemma} in the school of Landis, and it corresponds to the \emph{second lemma of De Giorgi} in the divergence-form theory. 

Moreover, by applying this ``second lemma'' (strong measure-to-pointwise estimate) at all scales and using a covering argument, one easily obtains the Harnack inequality, and this last part is abstract and independent of the equation itself.

Adapting these techniques to hypoelliptic equations of H\"ormander type (as in \cite{H}) has been a challenging open problem for a long time. It was, however, already known to experts that ingredients (B) and (C) could be extended in many scenarios, including kinetic equations. Therefore, the main obstacle to extending the argument has been the lack of an ABP inequality, or a suitable substitute, for equations of Kolmogorov type or equations in Carnot groups. In this note, we discuss a new idea obtained by OpenAI's GPT-6 Astra that provides a new way to derive ABP-like estimates and, consequently, a version of the Krylov–Safonov Harnack inequality for kinetic equations of Kolmogorov type. For further discussion of other ABP-like results, and for the original AI-generated material, see \cite{WeNeedABP}. The purpose of this note is to give a gentle exposition of the proof, identifying the new ideas. We explain it in the elliptic, parabolic and kinetic settings.

  \medskip

We consider the Kolmogorov equation with rough coefficients, namely,
\begin{equation}\label{Kinetic}
  Ku := \partial_t u + v \cdot \nabla_x u - \text{tr}(A(t,\,x,\,v)D_v^2u) = 0.
\end{equation}
Here $(t,\,x,\,v) \in (-1,\,0) \times B_1 \times B_1 \subset \mathbb{R}^{2n+1}$ and $A \in \text{Sym}_{n \times n}$ satisfies
\begin{equation}\label{UE}
  \lambda I \leq A \leq \Lambda I \quad \text{for some} \quad 0 < \lambda \leq \Lambda < \infty.
\end{equation}

The key missing ingredient in the kinetic setting was the following Aleksandrov estimate:
\begin{thm}\label{Aleksandrov}
  If $Ku = f \geq 0$ in $Q_1 = (-1,\,0) \times B_1 \times B_1$ and $u|_{\partial_{in}Q_1} = 0$, then
  \begin{equation}\label{Power}
    \sup_{Q_1} u \leq C\|f\|_{L^p(Q_1)}
  \end{equation}
  for some $C(n,\,\lambda,\,\Lambda) > 0$ and $p = p(n,\,\lambda,\,\Lambda) \in (1,\, \infty)$.
\end{thm}
Here $\partial_{in}$ denotes the incoming boundary.  In particular, $\partial_{in}Q_1=\{t=-1\}\cup\{|x|=1,\,v\cdot x<0\}\cup\{|v|=1\}$.

In the parabolic case, the corresponding estimate can be established with $p = n+1$. This is accomplished using convex analysis methods and the study of certain contact sets.
However, these methods do not seem to adapt easily to the hypoelliptic setting. The approach to Theorem \ref{Aleksandrov} is instead by duality. The estimate (\ref{Power}) is equivalent to a higher-integrability estimate for solutions to the adjoint equation. Namely, if
$$\begin{cases}
    \partial_t\rho + v \cdot \nabla_x\rho - \text{div}_v\text{div}_v(A(t,\,x,\,v)\rho) = f \geq 0 \text{ in } Q_1, \\
    \rho|_{\partial_{in}Q_1} = 0,
  \end{cases}$$
then
\begin{equation}\label{HigherIntegrability}
  \|\rho\|_{L^q(Q_1)} \leq C\|f\|_{L^1(Q_1)} \text{ for some } q(n,\,\lambda,\,\Lambda) > 1.
\end{equation}
The exponent $p$ in (\ref{Power}) is the H\"{o}lder conjugate of $q$. The adjoint viewpoint was also used in \cite{FS}.

\begin{rem}
  Technically, we wrote the backward (reflected in time and velocity) adjoint equation above. However, we will abuse terminology and simply refer to it as the adjoint equation in the sequel.
\end{rem}

For the purposes of proving the Harnack inequality, it actually suffices to prove the weaker estimate
\begin{equation}\label{LlogL}
  \int_{Q_1} \eta(\rho) \leq C(n,\,\lambda,\,\Lambda),
\end{equation}
under the additional assumption that $f$ has mass $1$, for some universal function $\eta$ satisfying
\begin{equation}\label{Superlinear}
  \lim_{s \rightarrow \infty} \frac{\eta(s)}{s} = \infty.
\end{equation}
We briefly explain why in the next section. Most of our exposition will focus on proving (\ref{LlogL}) for a choice $\eta(s)$ that is essentially the same as $s\log(s)$, which for reasons explained in Section \ref{Overview} is particularly convenient.

The proof uses divergence-form methods. The starting point is to study energy estimates for the adjoint equation. These contain a bad term involving a derivative of the coefficients, which at first seems insurmountable. However, it turns out that one can evolve in $(x,\,v)$ by an auxiliary flow in such a way that this bad term is the dissipation rate of a bounded quantity. After integrating in the auxiliary time, the term becomes harmless.

We will first illustrate how the proof works for elliptic and parabolic equations, which are simpler because the auxiliary flow (standard heat flow) commutes with the adjoint equation. We will then discuss the proof of (\ref{LlogL}) for kinetic equations, which is more complicated because one must account for new terms in the ``flowed'' adjoint equation coming from commutators with the drift. Finally, at the end we will discuss how, with a small adaptation of the method, one can obtain the stronger $L^q$ estimate (\ref{HigherIntegrability}). We conclude with a discussion of the exponent.

Let us also mention related literature on the regularity theory of equations in nondivergence form. Classical works on Krylov--Safonov regularity theory are \cite{MR790223,MR214905,MR126604,MR0147776,KS,Safonov1983Harnack}. The kinetic analogue of the Krylov--Safonov theorem was recorded as an open problem in \cite{MR4807233}. In \cite{MR3951695} a Harnack inequality under a Cordes–Landis condition or continuity assumption was proven, using the ideas developed in the elliptic and the parabolic case in Landis's book \cite{Landis}. The one-dimensional stationary problem was treated in \cite{zhu2026stationarykinetickrylovsafonovestimates}. For coefficients $A = A(t, v)$ in general dimensions, or $a = a(x, v)$ in one dimension, \cite{LMM} proves a kinetic Aleksandrov estimate. In \cite{Companion} it is also explained how to obtain Hölder continuity in the kinetic setting, and thus the entire kinetic Krylov--Safonov theory, from an ABP estimate.


\subsection*{Declaration of AI Use.}
The proof of the kinetic Krylov--Safonov theorem was found by OpenAI's GPT-6 Astra Ultra on September 6, 2026.
A link to the chat and the initial output is available here \cite{WeNeedABP}.
We do not claim credit for this result. Anyone could have obtained it by prompting the AI about it.
\begin{rem}
  A similar proof of the Aleksandrov estimate was made public in \cite{hao2026harnackestimatesnondivergencekinetic} on September 14, 2026.
  The precise involvement of AI in their proof was not made public.
\end{rem}

\subsection*{Acknowledgments}

C.~Mooney is supported by NSF Grant DMS-2143668. L.~Niebel is supported by SNSF Starting Grant TMSGI2\textunderscore226018. L.~Silvestre is supported by NSF grant DMS-2350263. 

\section{From \texorpdfstring{$L\log L$}{L log L} estimate to Harnack inequality}\label{LlogLHarnack}
Here we briefly explain why (\ref{LlogL})--(\ref{Superlinear}) suffice for proving the Harnack inequality. None of the key ideas are lost by restricting our discussion to the simpler case of elliptic equations. Below, we will always assume that $A$ satisfies (\ref{UE}).

In the proof of the Harnack inequality, the crucial estimate obtained using the Aleksandrov estimate is the following ``density estimate'':
$$\text{If} \quad \text{tr}(A(x)D^2u) \leq 0 \text{ in } B_1 \subset \mathbb{R}^n, \quad u \geq 0, \quad \text{and} \quad u(0) \leq 1,$$
then
\begin{equation}
  |\{u < C\}| \geq C^{-1} \text{ for some } C(n,\,\lambda,\,\Lambda) > 0.
\end{equation}

To see how this follows from (the elliptic analogue of) (\ref{LlogL})--(\ref{Superlinear}), note that by the maximum principle, $v := u + 2(|x|^2-1)$ is bounded from below by the solution $w$ to
$$\begin{cases}
    \text{tr}(A(x)D^2w) = 4n\Lambda\chi_{\{v < 0\}}, \\
    w|_{\partial B_1} = 0.
  \end{cases}$$
Moreover, $|w(0)| \geq 1$. Thus, if $\rho$ is the Dirichlet Green's function with pole at $0$, we have
$$1 \leq C(n,\,\Lambda)\int_{\{v < 0\}} \rho.$$
It is a simple exercise in real analysis to conclude from (\ref{LlogL}) and (\ref{Superlinear}) that
$$1 \leq C(n,\,\lambda,\,\Lambda)|\{v < 0\}| \leq C(n,\,\lambda,\,\Lambda)|\{u < 2\}|.$$

\begin{rem}
  Technically we are applying (\ref{LlogL}) with measure data, but this can be justified by approximation.
\end{rem}


\section{Overview}\label{Overview}

We now give an informal overview of the proof. We first do so in a schematic way, and then more technically in
the simplest setting of elliptic equations. We ignore technical issues like solution class, domains, etc.

\subsection{Main conceptual ideas}

The short story is that by working on the adjoint and along a well-chosen heat flow, one can reduce the Aleksandrov non-concentration estimate to the energy estimates of the divergence-form theory of De Giorgi--Nash--Moser: the derivatives falling on the rough coefficients in the diffusive directions are controlled by a ``two-solution'' entropy formula on the heat flow, and the derivatives in the transported directions are controlled by a choice of heat flow inspired by hypocoercivity, with vanishing coefficients for both diffusion in the transported variable and the cross term $\nabla_x \cdot \nabla_v$, in accordance with the scaling of the equation.

\smallskip
\textbf{Idea 1. Working with the adjoint.}

It is straightforward to see that establishing any universal non-concentration estimate for the adjoint fundamental solution that is strictly stronger than an $L^1$ bound yields an Aleksandrov estimate and allows us to apply the rest of the Krylov--Safonov--Tso machinery to obtain the Harnack inequality and Hölder continuity. Working with the adjoint then has the advantage of allowing us to commute the diffusive derivatives with regularizing semigroups and to combine the rough coefficients $A$ with the solution $\rho$ so that $A\rho$ can be treated as a new artificial unknown. The main (sub)-equation is then
\begin{equation*}
  \partial_t \rho + v \cdot \nabla_x \rho \le \nabla_v ^2 : (A \rho) + f.
\end{equation*}

\bigskip
\textbf{Idea 2. Reduce the non-divergence-form theory to the divergence-form theory}

The next idea is strange at first sight. It goes like this: ``if we had control over $\int \rho |\nabla \cdot A|^2$, then we could rewrite the equation in divergence form and obtain an energy estimate'':
\begin{equation*}
  \partial_t \rho + v \cdot \nabla_x \rho \le \nabla_v \cdot \left( A \nabla_v \rho \right) + \nabla_v \cdot \left( \rho \nabla_v \cdot A \right)+ f
\end{equation*}
because now the term $\nabla_v \cdot \left( \rho \nabla_v \cdot A \right)$ would be controlled by the standard Moser-type energy estimate, using the negative term created by the divergence term $\nabla_v \cdot \left( A \nabla_v \rho \right)$. The energy estimate here can be carried out in $L^p$ spaces or in terms of entropy, including more general convex entropies.

\smallskip
\textbf{Idea 3. Use a ``two-solution'' energy estimate on the heat flow}

Energy estimates for a single solution of the heat flow, formulated using any convex entropy, are well known, and we use them in the previous step as part of the energy estimate on our more general equation. We use, however, a more subtle energy estimate on the heat flow. This is the first of two really novel contributions. Consider two solutions $f_t$ and $g_t$ to the heat equation with $f_t>0$ on $\mathbb R^n$ or the torus $\mathbb T^n$. Then we have
\begin{equation*}
  \frac{{\rm d}}{{\rm d}t}\int \frac{g_t^2}{f_t} = - 2 \int f_t\left|\nabla \left( \frac{g_t}{f_t} \right) \right|^2.
\end{equation*}
In fact, the energy estimate between two unrelated solutions relies on an even stronger pointwise identity:
\begin{equation*}
  (\partial_t-\Delta)\left(\frac{g_t^2}{f_t}\right)
  =-2f_t\left|\nabla\left(\frac{g_t}{f_t}\right)\right|^2.
\end{equation*}
To see why this estimate is useful here, consider solutions with initial data $a\rho$ and $\rho$, where $a$ is rough. The right-hand side then controls, along the heat flow, some regularity of a quantity that is very close to $a$ for small $t$. This calculation can be found in \cite{plastino2} or \cite{MR2162224} and was used later in \cite{MR4748461}. The functional measures the Pearson \(\chi^2\)-divergence between the two evolving distributions.

\smallskip
\textbf{Idea 4. Regularise the equation along the heat flow.} Applying the heat flow and denoting its auxiliary time by $h$, we obtain the ``flowed equation''
\begin{equation*}
  \partial_t \rho_h + v \cdot \nabla_x \rho_h \le \nabla_v ^2 : (A \rho)_h + f_h.
\end{equation*}
Assume first that the heat flow only affects the diffusive directions $v$, and neglect the commutator error with the transport term. Note that the term is $(A \rho)_h$, rather than $A\rho_h$, because the heat flow does not commute with multiplication by $A$. Then define $B_h = (A\rho)_h /\rho_h$ and rewrite as (dropping the $h$ for now)
\begin{equation*}
  \partial_t \rho + v \cdot \nabla_x \rho \le \nabla_v \cdot \left( B \nabla_v \rho \right) + \nabla_v \cdot \left( \rho \nabla_v \cdot B \right)+ f.
\end{equation*}
The two-solution energy estimate controls $\int \rho |\nabla_v \cdot B|^2$  along $h$, which allows us to perform the standard Moser-like energy estimate. However, we still lack control in the transported variable $x$, and we have neglected the commutator error with the transport term.

\smallskip
\textbf{Idea 5. Hypocoercive modification of the heat flow.}

We now choose heat flow in both $x$ and $v$, with vanishing diffusion in $x$ in a way that respects the kinetic scaling as $h \to 0$. We also include \emph{cross-diffusion terms} in a way similar to the construction of a hypocoercivity norm \cite{MR2562709} with coefficients scaling as in \cite{MR2294477}: the coefficients are chosen to preserve the positivity of the quadratic form and to produce, through the commutator with the transport term, a negative $x$-diffusion term of higher order in $h$. This allows us to close the estimate.

\subsection{Concrete sketch in the elliptic case}

Assume that $A \in [\lambda I,\, \Lambda I]$ and that $\rho > 0$ solves
\begin{equation}\label{MotAdj}
  -\text{div div}(A(x)\rho) = f.
\end{equation}
Let $\eta$ and $\omega$ be functions on $[0,\,\infty)$ such that $\eta$ is convex and
\begin{equation}\label{etaomega}
  \eta''(s)s^2 \leq \omega(s),\, s > 0.
\end{equation}
If we test (\ref{MotAdj}) with $\eta'(\rho)$, integrate by parts, and use Cauchy--Schwarz, we get
\begin{equation}\label{MotEE}
  \int \eta''(\rho)|\nabla \rho|^2 \leq C(\lambda)\int \omega(\rho)|\text{div}A|^2 + C(\lambda)\int f\eta'(\rho).
\end{equation}
At first sight this estimate seems useless, since we have no control on the derivatives of $A$. The idea is to identify the left side and the first term on the right side as dissipation rates of entropy-like quantities along an auxiliary heat flow.

Let the subscript $h$ denote evolution by heat flow for time $h$. Flowing both sides of (\ref{MotAdj}) gives
\begin{equation}\label{MotFlowedAdj}
  -\text{div div}(B^h\rho_h) = f_h, \quad B^h := \frac{(A\rho)_h}{\rho_h} \in [\lambda I,\, \Lambda I].
\end{equation}
The idea of introducing this regularized diffusion coefficient goes back to \cite[Proof of Lemma 3.3]{Krylov1979PassingLimitII}.
Thus, by the same calculations leading to (\ref{MotEE}), we have
\begin{equation}\label{FlowedMotEE}
  \int \eta''(\rho_h)|\nabla \rho_h|^2 \leq C(n,\,\lambda)\int \omega(\rho_h)|\nabla B^h|^2 + C(\lambda)\int f_h\eta'(\rho_h).
\end{equation}
The left side is $-\frac{d}{dh}\int \eta(\rho_h)$, since
\begin{equation}\label{MotEnt1}
  (\partial_h - \Delta)(\eta(\rho_h)) = -\eta''(\rho_h)|\nabla \rho_h|^2.
\end{equation}
Recovering entropy by integrating its dissipation along an auxiliary flow/semigroup is a classical tool in the context of proving logarithmic Sobolev inequalities. We refer to \cite[Proposition 5 and its proof, Corollary 1]{MR889476}.

The main issue is to understand how the first term on the right arises as a dissipation. Note that the entries of $\rho_hB^h$ are caloric. The term $|\nabla B^h|^2$ will thus appear naturally in dissipations of quantities involving $|B^h|^2$. This leads us to the following question:

Assume that $\alpha$ and $\beta$ are functions such that $\alpha$ and $\alpha\beta$ are both caloric. For which $\omega$ can we guarantee that
\begin{equation}\label{MotEnt2}
  (\partial_h - \Delta)(\omega(\alpha)\beta^2) \leq -\kappa\omega(\alpha)|\nabla \beta|^2
\end{equation}
for some $\kappa > 0$?

A relatively short calculation reveals that the condition
\begin{equation}\label{omega}
  (2-\kappa)\omega\omega'' \geq 4(\omega'-\omega/s)^2,\, s > 0
\end{equation}
suffices. In this case, by applying (\ref{MotEnt2}) with $\alpha = \rho_h$ and $\beta$ any entry of $B^h$, we control the bad term by $-C(n,\,\lambda,\,\kappa)\frac{d}{dh}\int \omega(\rho_h)|B^h|^2$. Thus, after integrating (\ref{FlowedMotEE}) in $h$, we can hope to bound $\int \eta(\rho)$ provided
\begin{equation}\label{omegasmall}
  \omega(s) \text{ is small compared to } \eta(s), \text{ at least for large } s.
\end{equation}

There is some flexibility in choices of $\eta$ and $\omega$ satisfying (\ref{etaomega}), (\ref{omega}), (\ref{omegasmall}) for which the argument will close. The choice $$\eta(s) = s\log(s),\, \omega(s) = s$$
is convenient because the bad term is then the dissipation rate of a bounded quantity. We will use (a variant of) this choice to obtain the $L\log L$ estimate (\ref{LlogL}). Similarly, the choice
$$\eta(s) = s^q,\, \omega(s) = s^2\eta''(s)$$
works, provided $q > 1$ is sufficiently close to $1$. In a later section we will sketch how this choice gives the $L^q$ estimate (\ref{HigherIntegrability}).


\section{Entropies}\label{Entropies}
We now record a couple of entropy formulae that will be useful in the sequel, motivated by the discussion in the previous section. Assume that $\rho > 0,$ and that both $\rho$ and $\beta\rho$ satisfy the parabolic PDE
\begin{equation}\label{Aux}
  P_Mu := \partial_hu - \text{div}_y(M(h,\,y)\nabla_yu) = 0.
\end{equation}
Here $(h,\,y) \in \mathbb{R}_+ \times \mathbb{R}^d,\,M \in \text{Sym}_{d \times d},$ and $M \geq 0$.
Note that for any choice of $\eta$ we have
\begin{equation}\label{Entropy1}
  P_M(\eta(\rho)) =  -\eta''(\rho)|\nabla_y\rho|^2_M.
\end{equation}
Here and below, $|\xi|^2_M := M\xi\cdot\xi$. From here on out we fix
\begin{equation}
  \eta(s) := \begin{cases}
    s\log(s) - s + 1, \quad s > 1 \\
    0, \quad \text{otherwise}.
  \end{cases}
\end{equation}
Notice that $\eta = 0$ for $s \leq 1$. The corresponding choice of $\omega$ from the discussion in the previous section is $\omega(s) = s$. This clearly satisfies (\ref{etaomega}), (\ref{omega}) with $\kappa = 2$, and (\ref{omegasmall}), and moreover (\ref{MotEnt2}) becomes an equality:
\begin{equation}\label{Entropy2}
  P_M(\rho\beta^2) = -2\rho |\nabla_y\beta|^2_M.
\end{equation}


\section{Elliptic Case}\label{EllipticCase}
In this section and below, we will assume that all ingredients are smooth in order to emphasize ideas. The important thing is to obtain estimates that are independent of the regularity of the coefficients. The general case then follows by approximation.

Assume that $A(x) \in \text{Sym}_{n \times n}$ has eigenvalues in $[\lambda,\,\Lambda]$ and that $\rho$ solves
\begin{equation}\label{EllipticAdjoint}
  \begin{cases}
    -\text{div div}(A(x)\rho) = f \text{ in } B_1, \quad f \geq 0, \quad \int_{B_1} f = 1, \\
    \rho|_{\partial B_1} = 0.
  \end{cases}
\end{equation}
Then $\rho \geq 0$, by duality and the maximum principle. Extend $\rho,\,f$ to be zero outside $B_1$, and extend $A$ outside $B_1$ preserving (\ref{UE}). Using that $A\nabla\rho\cdot\nu \leq 0$ on $\partial B_1$ (here $\nu$ is the outer unit normal), we get
\begin{equation}\label{EllipticAdjointSub}
  -\text{div div}(A(x)\rho) \leq f \text{ in } \mathbb{R}^n.
\end{equation}
Let the subscript $h$ denote evolution by heat flow for time $h$. Then
$$B^h := \frac{(A\rho)_h}{\rho_h} \in [\lambda I,\, \Lambda I],$$
and we have for all $h > 0$ that
\begin{equation}\label{FlowedElliptic}
  -\text{div div}(B^h\rho_h) \leq f_h.
\end{equation}
It is easy to see that $\int_{B_1} \rho \leq C(n,\,\lambda)$ by testing (\ref{EllipticAdjoint}) with $1-|x|^2$, so we have the $L^1-L^{\infty}$ estimate
\begin{equation}\label{L1Linf}
  \|\rho_h\|_{L^{\infty}} \leq C(n,\,\lambda)h^{-n/2}.
\end{equation}
Testing (\ref{FlowedElliptic}) with $\eta'(\rho_h)$, integrating by parts, applying Cauchy--Schwarz, and using the fact that $f$ has mass $1$ gives
\begin{equation}\label{EllipticDissipation}
  \int \eta''(\rho_h)|\nabla \rho_h|^2 \leq C(n,\lambda)\int \rho_h|\nabla B^h|^2 + C(\lambda)\| \eta'(\rho_h)\|_{L^{\infty}}.
\end{equation}
Integrating (\ref{EllipticDissipation}) in $h$ until $\|\rho_h\|_{L^{\infty}} < 1$ and applying (\ref{Entropy1}) and (\ref{Entropy2}) with $d=n,\, M = I$, along with (\ref{L1Linf}), gives the desired estimate
$$\int_{B_1} \eta(\rho) < C(n,\,\lambda,\,\Lambda).$$


\section{Parabolic Case}\label{ParabolicCase}
The parabolic case is nearly identical to the elliptic one. We will evolve by an auxiliary heat flow, but on each fixed physical time slice $\{t = \text{const.}\}$.

Assume that $A(t,\,x) \in \text{Sym}_{n \times n}$ has eigenvalues in $[\lambda,\,\Lambda]$ and that $\rho$ solves
\begin{equation}\label{ParabolicAdjoint}
  \begin{cases}
    \partial_t\rho-\text{div}_x\text{div}_x(A(t,\,x)\rho) = f \text{ in } Q_1, \quad f \geq 0, \quad \int_{Q_1} f = 1, \\
    \rho|_{\partial_p Q_1} = 0.
  \end{cases}
\end{equation}
Here $Q_1 = (-1,\,0) \times B_1$ is a parabolic cylinder, and $\partial_p Q_1$ is the parabolic boundary.

Again, $\rho \geq 0$. Extending $\rho,\,f$ to be zero outside $Q_1$ and extending $A$ outside $Q_1$ while preserving (\ref{UE}), we get in the same way as above that
\begin{equation}\label{ParabolicAdjointSub}
  \partial_t \rho - \text{div}_x \text{div}_x(A(t,\,x)\rho) \leq f \text{ in } \mathcal{S} := (-1,\,0) \times \mathbb{R}^n.
\end{equation}
Let the subscript $h$ denote evolution by heat flow in $x$ for time $h$. For a physical time slice on which $\rho$ is not identically zero, let
$$B^h := \frac{(A\rho)_h}{\rho_h} \in [\lambda I, \,\Lambda I].$$
Otherwise, let $B^h = \lambda I$. Then
\begin{equation}\label{FlowedParabolic}
  \partial_t\rho_h - \text{div}_x \text{div}_x(B^h\rho_h) \leq f_h.
\end{equation}
It is easy to see that $\int_{B_1} \rho(t,\,x)\,dx \leq 1$ for all $t \in (-1,\,0)$ by integrating (\ref{ParabolicAdjointSub}), hence
$$\|\rho_h\|_{L^{\infty}} \leq C(n)h^{-n/2}.$$
Testing (\ref{FlowedParabolic}) with $\eta'(\rho_h)$, integrating by parts, and using Cauchy--Schwarz yields
\begin{equation}\label{ParabolicDissipation}
  \int \eta''(\rho_h)|\nabla \rho_h|^2 \leq C(n,\lambda)\int \rho_h|\nabla B^h|^2 + C(\lambda)\eta'(\|\rho_h\|_{L^{\infty}}).
\end{equation}
Here the integrals are over $\mathcal{S}$, and we have discarded a final-time integral with a favorable sign. The desired $L\log L$ estimate follows in the same way as in the elliptic case after integrating (\ref{ParabolicDissipation}) in $h$.


\section{Kinetic Case}\label{KineticCase}
Now assume that $A(t,\,x,\,v) \in \text{Sym}_{n \times n}$ has eigenvalues in $[\lambda,\,\Lambda]$ and that $\rho$ solves
\begin{equation}\label{KineticAdjoint}
  \begin{cases}
    \partial_t\rho + v \cdot \nabla_x\rho - \text{div}_v\text{div}_v(A\rho) = f \text{ in }Q_1, \quad f \geq 0, \quad \int_{Q_1} f = 1, \\
    \rho|_{\partial_{in} Q_1} = 0.
  \end{cases}
\end{equation}
Here $Q_1 = (-1,\,0) \times B_1 \times B_1$ is a unit kinetic cylinder, and
$$\partial_{in} Q_1 = \{t = -1\} \cup \{|x| = 1, v \cdot x < 0\} \cup \{|v| = 1\}.$$
Again, $\rho \geq 0$. Extend $\rho,\,f$ to be zero outside $Q_1$ and extend $A$ outside $Q_1$ while preserving (\ref{UE}). Arguing as above and using in addition that $\rho$ vanishes on the incoming spatial boundary (so the boundary contribution from the drift term has a favorable sign), we get
\begin{equation}\label{KineticAdjointSub}
  \partial_t \rho + v \cdot \nabla_x \rho - \text{div}_v \text{div}_v(A(t,\,x,\,v)\rho) \leq f \text{ in } \mathcal{S} := (-1,\,0) \times \mathbb{R}^n \times \mathbb{R}^n.
\end{equation}

We will evolve for time $h$ by an auxiliary flow in $x$ and $v$ of the type (\ref{Aux}), where $y = (x,\,v)$ and $d = 2n,$ for some carefully chosen $M(h,\,x,\,v)$. Velocity heat alone leaves position untouched, whereas adding position heat independently produces gradient terms that the equation does not control. The transport term provides a way to couple the two. A change in velocity changes the position reached by transport, so we allow a mixed position--velocity derivative in the smoothing. It is natural to choose a flow that respects the kinetic scaling $(h,\,x,\,v) \mapsto (r^2h,\,r^3x,\,rv)$.
For this reason, we consider diffusion matrices of the form
\begin{equation}\label{FlowChoice}
  M(h,\,x,\,v) = \begin{pmatrix} ah^2I & bhI \\ bhI & cI \end{pmatrix} := M^h, \quad a,\,b,\,c \in \mathbb{R}.
\end{equation}
Similar matrices appear in \cite{MR2294477} in the context of hypocoercivity and in \cite{MR3677826} in the context of a kinetic Bakry--\'Emery calculus.
We let the subscript $h$ denote evolution for time $h$ by the flow
$$P_M = \partial_h - \text{tr}(M^hD_{x,\,v}^2).$$
As above, we will let $B^h := (A\rho)_h/\rho_h$ on physical time slices for which $\rho$ is not identically zero, and define it to be $\lambda I$ otherwise.

Unlike in the elliptic and parabolic cases, new terms will appear in the ``flowed'' version of (\ref{KineticAdjointSub}) due to nontrivial commutators with the drift term:
\begin{equation}\label{FlowedKineticAdjoint}
  \partial_t\rho_h + v \cdot \nabla_x \rho_h - \text{div}_{x,\,v}(G^h\nabla_{x,v}\rho_h) - \text{div}_v(\rho_h\text{div}_vB^h) \leq f_h,
\end{equation}
where
\begin{equation}\label{FlowedDiffusion}
  G^h = \begin{pmatrix} -bh^2I & -chI \\ -chI & B^h \end{pmatrix}.
\end{equation}
We delay this calculation to the end of the section, and instead we immediately note that everything works as in the elliptic and parabolic cases provided the following key inequalities hold:
\begin{equation}\label{MhBound}
  ac > b^2, \quad a,\,c > 0
\end{equation}
and
\begin{equation}\label{GhMhBound}
  -\lambda b > c^2.
\end{equation}
The matrix $M^h$ determines the smoothing and the gradient terms in the entropy identities. The matrix $G^h$ determines what the smoothed equation can control. We therefore need $M^h$ to control velocity derivatives and $G^h$ to control $M^h$. Both requirements are satisfied by enforcing \eqref{MhBound} and \eqref{GhMhBound}.

Indeed, by integrating (\ref{KineticAdjointSub}) we see that $\int_{B_1 \times B_1} \rho(t,\,x,\,v)\,dx\,dv \leq 1$ for all $t \in (-1,\,0)$. The condition (\ref{MhBound}) thus implies the $L^1-L^{\infty}$ estimate
\begin{equation}\label{KL1Linf}
  \|\rho_h\|_{L^{\infty}} \leq C(a,\,b,\,c,\,n)h^{-2n},
\end{equation}
as well as decay in $x$ and $v$ of $\rho_h$. It also implies the ellipticity preservation $B^h \in [\lambda I,\, \Lambda I]$. These facts can be seen by using the fundamental solution of $P_M$, given by the Gaussian
\begin{equation}\label{KFundSoln}
  C(a,b,c,n)\,h^{-2n}\exp\left[-\frac{3c}{4ac-3b^2}\frac{\left|x-\frac{b}{2c}hv\right|^2}{h^3}-\frac{|v|^2}{4ch}\right].
\end{equation}

Moreover, it is not hard to see that (\ref{MhBound}) implies the existence of some $\kappa > 0$ such that
\begin{equation}\label{Mh2}
  M^h \geq \begin{pmatrix} 0 & 0 \\ 0 & \kappa I \end{pmatrix},
\end{equation}
and (\ref{GhMhBound}) implies in addition that
\begin{equation}\label{GhMh2}
  \quad G^h \geq \kappa M^h,
\end{equation}
up to taking $\kappa$ smaller if needed. Here we use that $B^h \geq \lambda I$.

Test (\ref{FlowedKineticAdjoint}) with $\eta'(\rho_h)$. What comes from the last term on the left side (which involves only $v$ derivatives) can be bounded above using (\ref{Mh2}), and the $|\nabla_{x,\,v} \rho_h|^2_{G^h}$ coming from the third term on the left side can be bounded below using (\ref{GhMh2}), giving
\begin{equation}\label{KineticDissipation}
  \int \eta''(\rho_h)|\nabla \rho_h|_{M^h}^2 \leq C\int \rho_h|\nabla B^h|_{M^h}^2 + C\eta'(\|\rho_h\|_{L^{\infty}}).
\end{equation}
The integrals are taken over $\mathcal{S}$, the constant $C = C(n,\,\lambda,\,a,\,b,\,c)$, and we have again discarded a final-time integral with a favorable sign. Integrating (\ref{KineticDissipation}) in $h$ and applying (\ref{Entropy1}), (\ref{Entropy2}), and (\ref{KL1Linf}) finishes the argument as before.

\begin{rem}
  There are many choices of $a,\,b,\,c$ satisfying (\ref{MhBound}) and (\ref{GhMhBound}), e.g. $a = \lambda,\, b = -\lambda/2,\, c = \lambda/2$ (which is similar to what Astra originally chose). More generally, take any $b < 0,\,c > 0$ satisfying (\ref{GhMhBound}), and then take $a > 0$ large enough that (\ref{MhBound}) is satisfied.
\end{rem}

\begin{rem}
  By the formula (\ref{KFundSoln}), we still have the $L^1-L^{\infty}$ bound (\ref{KL1Linf}) if $M^h$ is rank $1$ (when viewed as a $2 \times 2$ matrix), i.e. $ac = b^2$, provided $a,\,c > 0$ as well. This corresponds to auxiliary flow by a kinetic equation of Kolmogorov type, which is tempting as an initial guess. However, in this case, the inequality (\ref{Mh2}), which is vital for controlling derivatives of $B^h$ in the energy estimate, is not satisfied.
\end{rem}

To conclude the section we give a brief derivation of the flowed adjoint equation (\ref{FlowedKineticAdjoint}). Let $\mathcal{T} := \partial_t + v \cdot \nabla_x$. It suffices to prove, for any $g$, that
\begin{equation}\label{Commutator}
  (\mathcal{T} g)_h = \mathcal{T} g_h + bh^2\Delta_xg_h + 2ch\nabla_x \cdot \nabla_v g_h.
\end{equation}
The sides of (\ref{Commutator}) coincide at $h = 0$, so it suffices to show that the right side solves $P_M(\cdot) = 0$, i.e. that
$$[\partial_h - (ah^2\Delta_x + 2bh\nabla_x \cdot \nabla_v + c\Delta_v),\,\mathcal{T} + bh^2\Delta_x + 2ch\nabla_x \cdot \nabla_v] = 0.$$
This follows quickly from the commutator relations
$$[\Delta_v,\, \mathcal{T}] = 2\nabla_x \cdot \nabla_v, \quad [[\Delta_v,\, \mathcal{T}],\,\mathcal{T}] =[2\nabla_x \cdot \nabla_v,\, \mathcal{T}] = 2\Delta_x, \quad [\Delta_x,\,\mathcal{T}] = 0.$$
\begin{rem}
  When studying kinetic equations in non-divergence form we work with $\partial_t+v\cdot\nabla_x$ and $\Delta_v$ and spatial diffusion is recovered through two commutators.
  By contrast, for kinetic equations in divergence form we have access to $\nabla_v$ on its own and
  $$[\nabla_v,\,\partial_t+v\cdot\nabla_x]=\nabla_x$$
  already recovers the missing spatial directions.
\end{rem}

\begin{rem}
  The key coercivity term $-bh^2\Delta_x$ in the flowed equation (\ref{FlowedKineticAdjoint}) comes from the commutator of $\nabla_x \cdot \nabla_v$ with the drift. The bad mixed term $-2ch\nabla_x \cdot \nabla_v$ in (\ref{FlowedKineticAdjoint}) comes from the commutator of $\Delta_v$ with the drift. We are forced to take $c$ small to guarantee (\ref{GhMhBound}). However, the fact that $\Delta_x$ commutes with the drift allows us to choose $a$ freely to compensate for this and guarantee (\ref{MhBound}).
\end{rem}

\section{Aleksandrov Estimate}\label{AleksEstimate}
In this section, we indicate how to get from \eqref{LlogL} to the stronger $L^q$ estimate (\ref{HigherIntegrability}). The overall strategy is exactly the same as above, with slightly modified entropies and energy estimates. We sketch how it works in the elliptic setting; the parabolic and kinetic cases are entirely analogous.

  {\bf Entropies.} First, enforce the same assumptions on $\rho$ and $\beta$ as in Section \ref{Entropies}. We will update our choice of $\eta$ as follows:
$$\eta(s) = \begin{cases}s^q - qs + q - 1, \quad s > 1, \\
    0 \quad \text{ otherwise}.
  \end{cases}$$
Here $q > 1$. Motivated by the discussion in Section \ref{Overview}, we let
$$\omega(s) = q(q-1)s^q.$$
This choice clearly satisfies (\ref{etaomega}). It also satisfies (\ref{omega}) with constant
$$\kappa = \kappa_q = \frac{2(2-q)}{q},$$
for any $q \in (1,\,2)$. We record the associated entropy inequality, analogous to (\ref{MotEnt2}):
\begin{equation}\label{Entropy3}
  P_M(\omega(\rho)\beta^2) \leq -\kappa_q\omega(\rho)|\nabla_y\beta|^2_M.
\end{equation}
Finally, $\omega$ clearly satisfies (\ref{omegasmall}) provided $q$ is close to $1$. Below we will assume in particular that $q \in (1,\,3/2)$, which prevents $\kappa_q$ from degenerating (precisely, $\kappa_q > 2/3$).

  {\bf Elliptic Case.} Make the same assumptions on $A,\,\rho,\,f$ as in Section \ref{EllipticCase}. Test (\ref{FlowedElliptic}) with $\eta'(\rho_h)$. Integrating by parts and applying Cauchy--Schwarz gives
\begin{equation}\label{EE2}
  \int \eta''(\rho_h)|\nabla \rho_h|^2 \leq C(n,\,\lambda) \int \omega(\rho_h)|\nabla B^h|^2 + C(\lambda) \|\eta'(\rho_h)\|_{L^{\infty}}.
\end{equation}
Recalling the $L^1-L^{\infty}$ estimate (\ref{L1Linf}) and that $\eta(s) = 0$ for $s < 1$, we see that the integral from $0$ to $\infty$ in $h$ of the last term is bounded provided $q < 1 + 2/n$, which we now enforce. Performing this integration and recalling (\ref{Entropy1}) and (\ref{Entropy3}), we get
\begin{equation}\label{Almost}
  \int_{B_1} \eta(\rho) \leq C(n,\,\lambda,\,\Lambda)\int_{B_1} \omega(\rho) + C(n,\,\lambda,\,q).
\end{equation}
Finally, we recall that for $s \geq 0$,
$$s^q \leq \eta(s) + qs, \quad \omega(s) = q(q-1)s^q.$$
Combining these relations with (\ref{Almost}), we get
$$\int_{B_1} \rho^q \leq C(n,\,\lambda,\,\Lambda)q(q-1)\int_{B_1} \rho^q + q\int_{B_1}\rho + C(n,\,\lambda,\,q).$$
Taking $q(n,\,\lambda,\,\Lambda) < 1+2/n$ sufficiently close to $1$ we can absorb the first term on the right, and recalling the $L^1$ bound $\int_{B_1}\rho \leq C(n,\,\lambda)$, we get
$$\int_{B_1} \rho^q \leq C(n,\,\lambda,\,\Lambda).$$
For general $f$, the $L^q$ estimate (\ref{HigherIntegrability}) follows by homogeneity.

\begin{rem}
  There are two restrictions on $q$. The heat-kernel bound (\ref{KL1Linf}) controls the source integral only when $q<1+\frac{1}{2n}$ in the kinetic case. Absorbing the
  coefficient term imposes a second restriction, depending on ellipticity. The first reflects how the auxiliary flow spreads mass; the second comes
  from estimating the averaged coefficient by its size. It is this second restriction that makes our source exponent grow with the ellipticity ratio.
\end{rem}

\section{On exponents in the Aleksandrov estimate}
A careful tracking of constants in the argument of the previous section shows that we must take
$$q-1 \leq c(n)(\lambda/\Lambda)^2.$$
Thus, the exponent in the Aleksandrov estimate (\ref{Power}) produced by this argument goes to infinity with the ellipticity ratio $R := \Lambda/\lambda$.

On the other hand, in the elliptic and parabolic cases, the classical methods tell us much more about the Aleksandrov exponent.

For example, when $n \geq 2$, in the elliptic case we can take $p = n-\epsilon$ for some $\epsilon(n,\,R) > 0$ (see \cite{FS}). It was conjectured by Pucci \cite{MR208150} that the estimate holds for all
$$p > \max\left\{\frac{n}{2}, \frac{n}{1 + (n-1)/R}\right\}.$$
This threshold is predicted by radial homogeneous examples. It was proven to be sufficient by Astra \cite{WeNeedABP}.

In the parabolic case, we can take $p = n+1-\epsilon$ for some $\epsilon(n,\,R) > 0$ (\cite{FS}). Small modifications of the elliptic examples show that, for any exponent smaller than $n$, the Aleksandrov estimate fails at a sufficiently large ellipticity ratio.
Here, the condition
$$ p> \max\left\{\frac{n+2}{2}, \frac{n}{1 + (n-1)/R}\right\} $$
is necessary (see \cite[Theorem~4.8]{MR4861059}), while its sufficiency was established by Astra \cite{WeNeedABP}.

Finally, in the kinetic case, small modifications of stationary kinetically homogeneous examples show that, for any exponent smaller than $4n$, the Aleksandrov estimate fails at a sufficiently large ellipticity ratio. Moreover, under certain restrictions on the variables on which the coefficients depend, Aleksandrov estimates have been proven to hold for the best possible exponents.
These results, along with a more thorough discussion of optimal exponents and some explicit conjectures, can be found in \cite{LMM}.
The kinetic Pucci conjecture seems to be an interesting research direction as it presumably requires a novel alternative approach to the kinetic Aleksandrov estimate.

\raggedbottom
\bibliographystyle{abbrv}
\bibliography{ABP}
\vfill

\end{document}